\input amstex
\document
\magnification=\magstep1
\NoBlackBoxes 

\def\reel{\text{ }\hbox{{\rm R}\kern-.9em\hbox{{\rm I}}\text{ }}\text{ }}
\def\entier{\text{ }\hbox{{\rm N}\kern-.9em\hbox{{\rm I}}\text{ }}\text{ }}
\def\complexe{\text{ }\hbox{{\rm C}\kern-.55em\hbox{{\rm I}}\text{ }}}
\def\rationel{\text{ }\hbox{{\rm Q}\kern-.7em\hbox{{\rm I}}\text{ }}}
\def\un{\text{}\hbox{{\rm 1}\kern-.25em\hbox{{\rm I}}\,}}

\def\eps{\varepsilon}
\def\sn{\smallskip\noindent}
\def\mn{\medskip\noindent}
\def\bn{\bigskip\noindent}

\hbox{}
\centerline{\bf  LIPSCHITZ EMBEDDINGS OF METRIC SPACES INTO $c_0$}\medskip
\centerline{\bf  F. Baudier and R. Deville}

\bn\bf Abstract. 
\rm
let $M$ be a separable metric space.
We say that $f=(f_n):M\to c_0$ is a good-$\lambda$-embedding 
if, whenever $x,y\in M$, $x\ne y$ implies $d(x,y)\le\Vert f(x)-f(y)\Vert$ and,
for each $n$, $Lip(f_n)<\lambda$, where $Lip(f_n)$ denotes
the Lipschitz constant of $f_n$.
We prove that  there exists 
a good-$\lambda$-embedding from $M$ into $c_0$ if and only if 
$M$ satisfies an internal property called $\pi(\lambda)$.
As a consequence, we obtain that for any  separable metric space $M$, there exists 
a good-$2$-embedding from $M$ into $c_0$.
These statements slightly extend  former results obtained by N. Kalton and G. Lancien,
with simplified proofs.

\bn\mn
\bf 1) Introduction.

\mn\rm
First, let us recall that if $f$ is a mapping between the metric spaces $(M,d)$ and $(N,\delta)$,
the Lipschitz constant $Lip(f)$ is the infimum of all $\lambda$ such that for all $(x,y)\in M^2$, $\delta(f(x),f(y))\le\lambda d(x,y)$.

\sn
Let $(M,d)$ be a separable metric space and $\lambda\ge 1$. We say that 
$f:M\to c_0$ is a $\lambda$-embedding if, whenever $x,y\in M$, then~:
$$
d(x,y)\le\Vert f(x)-f(y)\Vert\le \lambda d(x,y).\tag 1
$$
Let us denote $f=(f_n)$ and, for each $n$, $E_n=\{(x,y)\in M\times M;\,d(x,y)\le\vert f_n(x)-f_n(y)\vert\}$.
Whenever $x,y\in M$, we have $\Vert f(x)-f(y)\Vert=\max\limits_n\vert f_n(x)-f_n(y)\vert$.
Hence $f$ is a $\lambda$-embedding if for each $n$, $Lip(f_n)\le\lambda$
and $M\times M=\bigcup\limits_nE_n$.

\mn
I. Aharoni [1] proved that for any separable metric space $M$, there exists a 
$\lambda$-embedding from $M$ into $c_0$ for any $\lambda>6$,
and that there is no $\lambda$-embedding from $\ell^1$ into $c_0$ if $\lambda<2$.
P. Assouad [2] improved this result by showing that one can construct a
$\lambda$-embedding from any separable metric space $M$ into $c_0$ for any $\lambda>3$.
Later on, J. Pelant [4] obtained the same result with $\lambda=3$.
It was also observed that there is no $\lambda$-embedding from $\ell^1$ 
into $c_0^+$ if $\lambda<3$.
All these authors actually constructed $\lambda$-embeddings 
into the positive cone $c_0^+$ of $c_0$.
Finally N. Kalton and G. Lancien [3] proved that for any separable metric space $M$, 
there exists a  $2$-embedding from $M$ into $c_0$, and this result is optimal 
(consider $M=\ell^1$).

\mn
We say that $f:M\to c_0$ is a strict-$\lambda$-embedding if, whenever $x,y\in M$
and $x\ne y$, then~:
$$
d(x,y)<\Vert f(x)-f(y)\Vert<\lambda d(x,y).\tag 2
$$



\mn
We say that $f=(f_n):M\to c_0$ is a good-$\lambda$-embedding if,
for each $n$, $Lip(f_n)<\lambda$, 
and $M\times M=\bigcup\limits_nE_n$.

\mn\bf 
Proposition 1.1. \it 
Assume $f:M\to c_0$ is a good-$\lambda$-embedding. Then,  
there exists $g:M\to c_0$ which is a strict and good-$\lambda$-embedding.

\mn
\sl Proof. \rm
Let $\lambda_n<\lambda$ be such that $f_n:M\to\reel$ is $\lambda_n$-Lipschitz continuous,
and let us define $g=(g_n):M\to c_0$ such that for each $n$, 
$g_n=\alpha_nf_n$ with $1<\alpha_n<2$ and $\alpha_n\lambda_n<\lambda$. 
Clearly, $g$ is still a good-$\lambda$-embedding.
If $x\ne y$, since the sequences $(f_n(x))$ and $(f_n(y))$ tend to zero, 
the sequence $\bigl(g_n(x)-g_n(y)\bigr)$ also converges to $0$ and there exists
$n_0$ such that
$$
\Vert g(x)-g(y)\Vert=\vert g_{n_0}(x)-g_{n_0}(y)\vert
\le\alpha_{n_0}\lambda_{n_0} d(x,y)<\lambda d(x,y)
$$ 
Since $\Vert f(x)-f(y)\Vert\le \Vert g(x)-g(y)\Vert$, this implies that 
$\Vert f(x)-f(y)\Vert<\lambda d(x,y)$. 
On the other hand, let $m_0$ be such that 
$\Vert f(x)-f(y)\Vert=\vert f_{m_0}(x)-f_{m_0}(y)\vert$.
We have
$$
d(x,y)\le\Vert f(x)-f(y)\Vert=\vert f_{m_0}(x)-f_{m_0}(y)\vert<
\vert\alpha_{m_0} f_{m_0}(x)-\alpha_{m_0}f_{m_0}(y)\vert\le\Vert g(x)-g(y)\Vert.
$$
Therefore, $g$ is also a strict-$\lambda$-embedding.

\mn
Our purpose is to prove that for every separable metric space, one can construct a strict-2-embedding from 
$M$ into $c_0$. We introduce also  a property $\pi(\lambda)$ of a metric space,
slightly weaker than a property  introduced by N. Kalton and G. Lancien, 
and we prove that if $1<\lambda\le 2$, a separable metric space $M$ admits a good-$\lambda$-embedding into $c_0$ if and only if it has the property $\pi(\lambda)$.

\bn\mn\bf
2) Necessary condition for the existence of good-$\lambda$-embedding into $c_0$.

\mn\rm
Let $(M,d)$ be a metric space and $E$ be a non empty subset of  $M\times M$.
We denote $\pi_1(E)=\{x\in M;\,\exists y\in M, (x,y)\in E\}$, $\pi_2(E)=\{y\in M;\,\exists x\in M, (x,y)\in E\}$
the projections of $E$, and
$\pi(E)=\pi_1(E)\times  
\pi_2(E)$ the smallest rectangle containing $E$.
We also define the gap of $E$ by $\delta(E):=\inf\{d(x,y);(x,y)\in E\}$
and the diameter of $E$ by $diam(E)=\sup\{d(x,y);\,(x,y)\in E\}$. 
These notions are not quite standard, and require some comments. 
Let us denote $\Delta:=\{(x,x);\,x\in M\}$ the diagonal of $M\times M$, 
and let us endow the set $M\times M$ with the metric 
$d_1\bigl((x,y),(x',y')\bigr)=d(x,x')+d(y,y')$.
The distance from a point $(y,z)\in M\times M$ to $\Delta$
is 
$$
d_1\bigl((y,z),\Delta\bigr)=\inf\{d_1\bigl((y,z),(x,x)\bigr);\,(x,x)\in\Delta\}
$$
and it is easy to check that $d_1\bigl((y,z),\Delta\bigr)=d(y,z)$.
Consequently, if $\emptyset\ne E\subset M\times M$, the 
smallest distance from a point of $E$ to $\Delta$
is the quantity 
$$
d_1(E,\Delta)=\inf\{d_1\bigl((y,z),\Delta\bigr);\,(y,z)\in E\}=\delta(E)
$$
On the other hand, the largest distance from a point of $E$ to $\Delta$ is 
$$
D_1(E,\Delta)=\sup\{d_1\bigl((y,z),\Delta\bigr);\,(y,z)\in E\}=diam(E)
$$
Whenever $E$ is of the form $U\times V$, then $\delta(E)=\inf\{d(x,y);\,x\in U,\,y\in V\}$
is the gap between $U$ and $V$,
and $diam(E)=\sup\{d(x,y);\,x\in U,\,y\in V\}$.
Thus, if $U=V$, $diam(E)$ is the usual diameter of $U$.

\bn\bf Fact 2.1. \it Let  $E$ be a bounded subset of $M\times M$, $F$ be a finite dimensional normed vector space, 
let $P:M\to F$ be
such that $Lip(P)\le\lambda$ and $d(x,y)\le\Vert P(x)-P(y)\Vert$ for each $(x,y)\in E$, and let $\eps>0$.  
Then, there exists a finite partition $\{E_1,\cdots, E_N\}$ 
of $E$ so that 
$$
\text{ for each }n ,\quad diam(E_n)<\lambda\delta(\pi(E_n))+\eps.
$$

\mn\sl Proof. \rm 
The set $P(\pi_1(E)\cup\pi_2(E))\subset F$ is bounded as $E$ is bounded 
and $P$, $\pi_1$ and $\pi_2$ are Lipschitz.
Hence we can find a finite partition of this set into subsets $F_j$
of diameter $<\eps/4$. The sets
$E_{j,k}=(P^{-1}(F_j)\times P^{-1}(F_k))\cap E$ which are non empty form a partition of  $E$.
If $(x,y)\in E_{j,k}$ and $(u,v)\in\pi(E_{j,k})$, then
$$
\Vert P(x)-P(y)\Vert\le\Vert P(x)-P(u)\Vert+\Vert P(u)-P(v)\Vert+\Vert P(v)-P(y)\Vert\le\eps/2+\lambda d(u,v)
$$
Thus $d(x,y)\le\lambda d(u,v)+\eps/2$. 
The result follows by taking the infimum over all  $(u,v)\in\pi(E_{j,k})$,
the supremum over all $(x,y)\in E_{j,k}$, and by relabeling the sets $E_{j,k}$.

\bn\bf
Definition 2.2. \rm A metric space $(M,d)$ has property $\pi(\lambda)$
if, for any balls $B_1$ and $B_2$
of radii $r_1$ and $r_2$ and for any  non empty subset  $E$  of $B_1\times B_2$
satisfying $\delta(E)>\lambda(r_1+r_2)$,  
there exists a partition $\{E_1,\cdots,E_N\}$ of $E$, 
such that
$$
\text{ for each }n ,\quad
diam(E_n)< \lambda\delta(\pi(E_n))
$$
We say that $(M,d)$ has the property weak-$\pi(\lambda)$ if  the conclusion is replaced by the weaker conclusion :
there exists non empty closed subsets $F_1,\cdots, F_N$ 
covering $E$ such that
$$
\text{ for each }n ,\quad
r_1+r_2<\delta(\pi(F_n))
$$
this conclusion is indeed weaker : if $F_n$ is the closure of $E_n$ in $M\times M$, then 
$\lambda(r_1+r_2)\!<\!\delta(E)\!\le\! diam(E_n)\!<\! \lambda\delta(\pi(E_n))\!=\!\lambda\delta(\pi(F_n))$.
It is also easy to see that if $\lambda<\mu$ and if $M$ has $\pi(\lambda)$,
then $M$ has $\pi(\mu)$,
and that if $M$ has at least $2$ elements, $M$ never has $\pi(1)$.

\bn\bf
Proposition 2.3. \it 1) Assume that there is a good-$\lambda$-embedding from $M$ into $c_0$.
Then $M$ has property~$\pi(\lambda)$.

\sn
2) If $(M,d)$ $\lambda$-embeds into $c_0$, then
$M$ has property weak-$\pi(\lambda)$.

\mn\sl Proof. \rm
Let $f:M\to c_0$ be a $\lambda$-embedding. If $(e_i)$ is the unit vector basis of $c_0$,
then $f(x)=\sum\limits_{i=0}^{+\infty}f_i(x)e_i$.  
Let $B_1$ and $B_2$ be balls of radii $r_1$ and $r_2$ and of centers $a_1$ and $a_2$, 
and $E\subset B_1\times B_2$ such that $\delta(E)>\lambda(r_1+r_2)$. We claim that  
the function $E\ni(x,y)\mapsto\Vert f(x)-f(y)\Vert$ depends on
finitely many coordinates, i. e. there exists $i_0\in\entier$ such that,
if $P(x)=\sum\limits_{i=0}^{i_0}f_i(x)e_i$ then $\Vert f(x)-f(y)\Vert=\Vert P(x)-P(y)\Vert$.

\noindent
Fix $\eps>0$ such that $\eps<\delta(E)-\lambda(r_1+r_2)$. We choose
$i_0$ such that, if $Q=f-P$, then $\Vert Q(a_1)-Q(a_2)\Vert<\eps$.
If $(x,y)\in E$, then
$$
\eqalign{\Vert Q(x)-Q(y)\Vert&\le\Vert Q(x)-Q(a_1)\Vert+
\Vert Q(a_1)-Q(a_2)\Vert+\Vert Q(a_2)-Q(y)\Vert\cr
&<\lambda(r_1+r_2)+\eps<\delta(E)\le d(x,y).
\cr}
$$
Hence
$d(x,y)\le\Vert f(x)\!-\!f(y)\Vert=\max\{\Vert Q(x)-Q(y)\Vert,\Vert P(x)\!-\!P(y)\Vert\}=\Vert P(x)\!-\!P(y)\Vert$.
This proves our claim.
Since $Lip(P)\le\lambda$, Fact 2.1 implies  
the existence of a a partition $\{E_1,\cdots, E_N\}$ 
of $E$ such that
for all $n$,\quad
$diam(E_n)<\lambda\delta(\pi(E_n))+\eps$.
Since we also have
$\lambda(r_1+r_2)+\eps\le diam(E_n)$,
we have $r_1+r_2<\delta(\pi(E_n))$, so if $F_n$ is the norm closure of $E_n$ in $E$,
$$
r_1+r_2<\delta(\pi(F_n))
$$
When $f$ is a good-$\lambda$-embedding, the mapping $P$
is $\mu$-Lipschitz continuous for some $\mu<\lambda$, so we can asume that for all $n$,
$diam(E_n)<\mu\delta(\pi(E_n))+\alpha$, where $\alpha=\min\{\eps,(\lambda-\mu)(r_1+r_2)\}$.
This still implies $r_1+r_2<\delta(\pi(E_n))$.
Finally,
$$
diam(E_n)<\mu\delta(\pi(E_n))+(\lambda-\mu)(r_1+r_2)<\lambda\delta(\pi(E_n)).
$$

\mn\bf
Corollary 2.4 (see [3]). \it Let $X$ be a Banach space. If there exists $u\in S_X$
and an infinite dimensional subspace $Y$ of $X$ such that
$\inf\{\Vert u+y\Vert;\,y\in S_Y\}>\lambda$,
then there is no $\lambda$-embedding from $X$ into $c_0$.

\mn\sl Proof. \rm
If $E=\{(u+y,-u-y);\,y\in S_Y\}\subset\overline{B}(u,1)\times\overline{B}(-u,1)$, 
$E$ satisfies $\delta(E)>2\lambda$.
Assume there exists a  $\lambda$-embedding from $M$ into $c_0$. 
Then $X$ has the weak-$\pi(\lambda)$ property, so there exists closed subsets
$F_1,\cdots,F_N$ of $E$  covering $E$ such that for each $n$,\quad
$\delta(\pi(F_n))>2$.

\sn
On the other hand, 
$A_n=\{y\in S_Y;\,(u+y,-u-y)\in F_n\}$ is closed
and $A_1\cup\cdots\cup A_N=S_Y$.
Since $dim(Y)>N$,  the Borsuk-Ulam thoeorem yields the existence of  
$y\in S_Y$ and $n$ such that $\{y,-y\}\subset A_n$.
Hence $(u+y,-u+y)\in\pi(F_n)$ and so $\delta(\pi(F_n))\le 2$,
which is absurd.

\bn\bf Example 2.5. \it There is no $\lambda$-embedding from $\ell^p$
into $c_0$ for any $\lambda<2^{1/p}$. In particular, $\ell^1$ is a metric space 
which does not  $\lambda$-embed into $c_0$ with 
$\lambda<2$. \rm (If $u=e_0$ and $Y=\{y=(y_i)\in\ell^p;\,y_0=0\}$, then 
$\Vert u+y\Vert=2^{1/p}$ for all $y\in S_Y$).

\bn\mn\bf
3) Examples of spaces with property $\pi(\lambda)$.\rm

\bn\bf
Example 3.1. \it  A metric  space such that the bounded subsets of $M$ are totally bounded has property
$\pi(1+\eps)$ for all $\eps>0$ \rm (partition $E$ into subsets $E_n$
of small $d_1$-diameter).

\bn\bf
Example 3.2. \it If $(M,d)$ is a metric space, 
then $(M,d)$ has property $\pi(2)$.

\noindent
Therefore, property $\pi(\lambda)$ is of interest only if $1<\lambda\le 2$.

\mn\sl Proof. \rm
Let   $E\subset B_1\!\times\! B_2$, with $B_1$ et $B_2$ balls 
of radii $r_1\ge r_2$, 
and assume $\eps:=\delta(E)-2(r_1+r_2)>0$.
Let $a_0=\delta(E)\!<\!a_1<\cdots<a_{N-1}\!<\!diam(E)\!<\!a_N$ so that 
for all $1\le n\le N$,  $a_n-a_{n-1}<\eps$.
Define, for $1\le n\le N$, 
$E_n=\big\{(x,y)\in E;\, a_{n-1}\le d(x,y)< a_{n}\big\}$.
Thus $\delta(E_n)+\eps>diam(E_{n})$.
If $(u,v)\in \pi(E_n)$, one can find $v'\in B_2$
such that $(u,v')\in E_n$. Moreover $v,v'\in B_2$, hence :  
$$
diam(E_n)<2\delta(E_n)-2(r_1+r_2)\le 2d(u,v')-2d(v',v)\le 2d(u,v)
$$
Taking the infimum over all $(u,v)\in \pi(E_n)$, we get
$diam(E_{n})<2\delta(\pi(E_n))$.

\bn\bf
Example 3.3. \it If $(X_n)$ is a sequence of finite dimensional Banach spaces, then $(\oplus X_n)_p$ has property $\pi(2^{1/p})$.

\mn\rm
\sl Proof. \rm Let   
$E\subset B(a_1,r_1)\times B(a_2,r_2)$ such that $\alpha=\delta(E)^p-2(r_1+r_2)^p>0$. 
We select
$\eps>0$ so that $(r_1+r_2+\eps)^p<\delta(E)^p/2-\alpha/4$ and 
$2(t+\eps)^p-\alpha/2<2t^p$ whenever $0\le t\le diam(E)$. 
If $x\in(\oplus X_n)_p$, then $x=\sum\limits_{n=1}^\infty x_n$ with $x_n\in X_n$ for each $n$.
Define $P,Q:(\oplus X_n)_p\to(\oplus X_n)_p$
by $P(\sum\limits_{i=0}^\infty x_i)\!=\!\sum\limits_{i=0}^{i_0} x_i$
and $Q\!=\!I\!-\!P$,
where $i_0$ is such that $\Vert Qa_1- Qa_2\Vert<\eps$. 
According to Fact 2.1, since $P$ is an operator of norm $1$ with values in a finite dimensional subspace of $(\oplus X_n)_p$, 
we can find relatively closed subsets $E_n$ of $E$ covering $E$ such that, 
for all $n$, if $(x,y)\in E_n$, then
$\Vert Px-Py\Vert\le\delta(\pi(E_n))+\eps$. 
On the other hand, 
$\Vert Qx -Qy\Vert\le\Vert Qx-Qa_1\Vert+\Vert Qy-Qa_2\Vert+\Vert Qa_1- Qa_2\Vert\le r_1+r_2+\eps$.
Moreover,
$\Vert x-y\Vert^p=\Vert Px-Py\Vert^p+\Vert Qx-Qy\Vert^p$, so
$$
diam(E_n)^p\le\bigl(\delta(\pi(E_n))+\eps\bigr) ^p+(r_1+r_2+\eps)^p<\bigl(\delta(\pi(E_n))+\eps\bigr) ^p+diam(E_n)^p/2-\alpha/4
$$
which implies $diam(E_n)^p<2\bigl(\delta(\pi(E_n))+\eps\bigr) ^p-\alpha/2<2\delta(\pi(E_n))^p$.

\bn\bf
Remark 3.4. \rm
In the definition of property $\pi(\lambda)$, if $a_1$ and $a_2$ are  centers of 
$B_1$ and $B_2$, we can assume that 
$$
(\lambda-1)(r_1+r_2)< d(a_1,a_2)\le\frac{\lambda+1}{\lambda-1}(r_1+r_2)
$$
Indeed, if $E\ne\emptyset$, then, for each $(x,y)\in E$, 
$\lambda(r_1+r_2)<\delta(E)\le d(x,y)\le d(a_1,a_2)+(r_1+r_2)$,
which proves the first inequality. If $\frac{\lambda+1}{\lambda-1}(r_1+r_2)<d(a_1,a_2)$,
the conclusion of property $\pi(\lambda)$ is always true if we take $N=1$ and $E_1=E$, since then
$$
diam(E)\le d(a_1,a_2)+(r_1+r_2)<\lambda\bigl(d(a_1,a_2)-(r_1+r_2)\bigr)\le\lambda\delta(\pi(E))
$$

\bn\mn\bf
4) Constructing strict-$\lambda$-embeddings into $c_0$.

\bn\bf
Theorem 4.1. \it Let $(M,d)$ be a separable metric space and $1\!<\!\lambda\!\le\! 2$. If $M$ has property  $\pi(\lambda)$,
there exists a $\lambda$-embedding from $M$ into $c_0$ which is strict and good. 

\mn\rm
Theorem 4.1 and Proposition 2.3 show that there exists a good-$\lambda$-embedding from $M$ into $c_0$ if and only if $M$ has property  $\pi(\lambda)$. We do not know any internal characterization of separable metric spaces 
that admit a $\lambda$-embedding into $c_0$.

\bn\bf
Corollary 4.2. \it Let $(M,d)$ be a metric space such that the bounded subsets of $M$ are totally bounded. 
For all $\eps>0$, there exists a $(1+\eps)$-embedding from $M$ into $c_0$.

\bn\bf
Corollary 4.3. \it Every separable metric space strictly-$2$-embeds into $c_0$.

\bn\rm
This result is optimal since $\ell^1$ does not
$\lambda$-embed into $c_0$ whenever 
$\lambda<2$.

\bn\bf
Corollary 4.4. \it   If $(X_n)$ is a sequence of finite dimensional Banach spaces, then $(\oplus X_n)_p$ strictly-$2^{1/p}$-embeds into $c_0$.

\bn\rm
This result is optimal since we have seen that there is no 
$\lambda$-embedding from $\ell^p$
into $c_0$ with $\lambda<2^{1/p}$.
We now turn to the proof of Theorem 4;1. We need some further notations.
If $(M,d)$ is a metric space , $x\in M$ and $U,V\subset M$, 
the distance from $x$ to $U$ is
$d(x,U)=\delta(\{x\}\times U)$ and 
the gap between $U$ and $V$ is $\delta(U,V)=\delta(U\times V)$.
The  coordinates  of the embedding from $M$ into $c_0$
are of the following type~:

\bn\bf
Lemma 4.5. \it Let $(M,d)$ be a metric space, $U,V,F$ three non empty subsets of $M$ 
and $\eps\ge 0$.
There exists $f:M\to\reel$, $1$-Lipschitz, such that~:

\medskip
1) For all $x\in F$, \quad $\vert f(x)\vert\le\eps$.

\medskip
2) For all $(x,y)\in U\times V$, \quad $f(x)-f(y)=\min\big\{\delta(U,V),\delta(U,F)+\delta(V,F)+2\eps\big\}$.

\bn\bf
Lemma 4.6. \it Let $1<\lambda\le 2$, $(M,d)$ be a metric space with property $\pi(\lambda)$, 
$F\subset G$ be finite subsets of $M$ and $0<\alpha<\beta$. We set~:
$$
A(F,\beta)=\big\{(x,y)\in M\times M;\,
\lambda\bigl(d(x,F)+d(y,F)+\beta\bigr)\le d(x,y)\big\}
$$
Then there exists a finite partition $\{E_1,\cdots,E_N\}$ of $A(G,\alpha)\backslash A(F,\beta)$
such that, if we denote $\pi(E_n)=U_n\times V_n$
then
$$
\text{for each }n\qquad diam(E_n)<\lambda\min\big\{\delta(U_n,V_n),\delta(U_{n},F)+\delta(V_{n},F)+2\beta\big\}.
$$

\mn\sl
Proof of Theorem 4.1.

\mn\rm
The goal is to construct a sequence $(f_n)$ of  $1$-Lipschitz continuous functions
satisfying, for every $x\in M$, $\lim\limits_{n\to\infty}f_n(x)=0$,
and a partition $\{E_n;\,n\in\entier\}$ of  $\{(x,y)\in M\times M;\,x\ne y\}$,
so that for each $n$, 
the function $(x,y)\to f_n(x)-f_n(y)$ is equal to some constant $c_n$ on $E_n$ and $diam(E_n)<\lambda c_n$. 
The required strict and good-$\lambda$-embedding is then $f=(\lambda_nf_n)$
where $\lambda_n<\lambda$ is chosen so that $diam(E_n)<\lambda_n c_n$.

\mn
Let $(a_k)$ be a dense sequence of distinct points of $M$, $F_k=\{a_1,\cdots,a_k\}$,
and $(\eps_k)$ be a decreasing sequence of real numbers converging to $0$.
We set 
$\Delta_k=A(F_{k+1},\eps_{k+1})\backslash A(F_k,\eps_k)$. 
The sets $\Delta_k$ form a partition of $\{(x,y)\in M\times M;\,x\ne y\}$.
Indeed, if $x,y\in M$, $x\ne y$ and if 
$\sigma_k=\lambda\bigl(d(x,F_k)+d(y,F_k)+\eps_k\bigr)$, then 
$0<d(x,y)<\sigma_1$, $(\sigma_k)$ is decreasing and $\lim\limits_{k\to\infty}\sigma_k=0$, so there exists a unique $k$ such that 
$\sigma_{k+1}\le d(x,y)<\sigma_k$, which means $(x,y)\in \Delta_k$.

\sn
By Lemma~4.6, there exists integers $0=n_1<n_2<\cdots<n_k<\cdots$ 
and subsets $E_n$ of $M\times M$ such that for~all~$k$,
$\{E_n;\,n_k< n\le n_{k+1}\}$ is a partition of $\Delta_k$, and, 
whenever $n_k< n\le n_{k+1}$ then
$$
diam(E_n)<\lambda\min\big\{\delta(U_n,V_n),\delta(U_n,F_k)+\delta(V_n,F_k)+2\eps_k\big\}.
$$
where $\pi(E_n)=U_n\times V_n$.  In particular,
$\{E_n;\,n\in\entier\}$ is a partition of  $\{(x,y)\in M\times M;\,x\ne y\}$.
By Lemma 4.5, there are  $1$-Lipschitz functions
$f_n:M\to\reel$ so that 

1) if $x\in F_k$ and $n_k< n\le n_{k+1}$, then $\vert f_n(x)\vert\le\eps_k$,

2) if $n_k< n\le n_{k+1}$ and $(x,y)\in U_n\times V_n$, then 
$$
f_n(x)-f_n(y)=c_n:=\min\big\{\delta(U_n,V_n),\delta(U_n,F_k)+\delta(V_n,F_k)+2\eps_k\big\}.
$$
and so $diam(E_n)<\lambda c_n$.

\mn
If $x\in M$, let us show that $\lim\limits_{n\to\infty}f_n(x)=0$.
If $\eps>0$, we fix $j$ such that $d(x,a_j)<\eps/2$,
then $k\ge j$ such that $\eps_k<\eps/2$. 
Since the functions $f_n$ are $1$-Lipschitz continuous, 
if $n\ge n_k$, 
then $\vert f_n(x)\vert\le d(x,a_j)+\vert f_n(a_j)\vert<\eps/2+\eps_k<\eps$.

\bn\sl
Proof of Lemma 4.5. \rm 
We fix $s,t$ such that  $-\delta(V,F)-\eps\le s\le 0\le t \le \delta(U,F)+\eps$ and 
$t-s=\min\big\{\delta(U,V),\delta(U,F)+\delta(V,F)+2\eps\big\}$, and we set
$$
f(x):=\min\big\{d(x,U)+t,d(x,V)+s,d(x,F)+\eps\big\}
$$
The function $f$ is $1$-Lipschitz continuous as the infimum of 
$1$-Lipschitz continuous functions.

\mn
If $x\in U$, $f(x)=\min\big\{t,d(x,V)+s,d(x,F)+\eps\big\}=t$ because
$d(x,V)+s\ge \delta(U,V)+s\ge t$ and $d(x,F)+\eps\ge\delta(U,F)+\eps\ge t$. 
If $y\in V$, $f(y)=\min\big\{d(y,U)+t,s,d(x,F)+\eps\big\}= s$ because $s\le 0$.
Therefore, if $x\in U$ and $y\in V$,
then $f(x)-f(y)=t-s$, which proves 2).

\mn
Finally, if $x\in F$, then $f(x)=\min\big\{d(x,U)+t,d(x,V)+s,\eps\big\}\le\eps$.
On the other hand $d(x,U)+t\ge 0$ et $d(x,V)+s\ge\delta(V,F)+s\ge-\eps$, 
so $f(x)\ge -\eps$, which proves 1).

\bn\sl
Proof of Lemma 4.6.

\mn\rm
Set $\Delta:=A(G,\alpha)\backslash A(F,\beta)$. 
There is a bounded subset $B$ of $M$ such that $\Delta\subset B\times B$, 
because $\lambda>1$, $G$ is bounded and   
$$
\lambda\bigl(d(x,G)+d(y,G)\bigr)\le d(x,y)\le d(x,G)+d(y,G)+diam(G).
$$
whenever $(x,y)\in A(G,\alpha)$. 
Thus, there is a partition $\{B_1,B_2,\cdots,B_m\}$ of the bounded set $B$  
such that for all $j$, if $x,x'\in B_j$ and $a\in G$, then  
$\vert d(x,a)-d(x',a)\vert\le\alpha/5$,
and so, 
$$
\text{for all }x\in B_j,\quad\text{for all }a\in G,  \qquad d(x,a)< d(B_j,a)+\alpha/4.
$$
Since $G$ is finite, there exists $a_j\in G$ such that $d(B_j,a_j)=\delta(B_j,G)$, 
and so
$B_j\subset B(a_j,r_j)$, where  $r_j=\delta(B_j,G)+\alpha/4$.
The subsets $E_{jk}=\Delta\cap B_j\times B_k$ of $\Delta$ 
form a partition of $\Delta$,
$E_{jk}\subset B(a_j,r_j)\times B(a_k,r_k)$, and, if $(x,y)\in\ E_{jk}$~:
$$
\eqalign{
d(x,y)&\ge \lambda\bigl(d(x,G)+d(y,G)+\alpha\bigr)\ge \lambda\bigl(\delta(B_j,G)+\delta(B_k,G)+\alpha\bigr)\cr
&= \lambda(r_j+r_k+\alpha/2)> \lambda(r_j+r_k).\cr
}
$$
So $\delta(E_{jk})> \lambda(r_j+r_k)$. 
According to property $\pi(\lambda)$ applied to each $E_{jk}$, 
there exists  a finite partition
$\{E_1,\cdots,E_N\}$ of $\Delta$ such that, 
$$
diam(E_n)< \lambda\delta(\pi(E_n))=\lambda\delta(U_n,V_n),
$$
where $\pi(E_n)=U_n\times V_n$. Moreover,
if $j,k,n$ are such that $E_n\subset B_j\times B_k$ and if $(x,y)\in E_n$, then
$$
\eqalign{
d(x,y)\le\lambda\bigl(d(x,F)+d(y,F)+\beta\bigr)&
\le\lambda\bigl(\delta(B_j,F)+\alpha/4+\delta(B_k,F)+\alpha/4+\beta\bigr)\cr&
\le\lambda\bigl(\delta(U_n,F)+\delta(V_n,F)+\alpha/2+\beta\bigr).
\cr}
$$
hence
$$
diam(E_n)\le \lambda\bigl(\delta(U_n,F)+\delta(V_n,F)+\alpha/2+\beta\bigr)<\lambda\bigl(\delta(U_n,F)+\delta(V_n,F)+2\beta\bigr).
$$

\bn\bf
5) Some consequences.

\bn\rm
Observe that a metric space has property $\pi(\lambda)$ (resp. weak-$\pi(\lambda)$)
if and only if its bounded subsets have it. In particular a Banach space has property 
$\pi(\lambda)$ (resp. weak-$\pi(\lambda)$)
if and only if its unit ball has it.
Since the property ``there exists a good-$\lambda$-embedding from $M$ into $c_0$''
is equivant to the property 
``$M$ has $\pi(\lambda)$'', we can state~:

\bn\bf
Proposition 5.1. \it Assume that $(M,d)$ is a separable metric space and that for each ball $B$ of $M$,
there is a good-$\lambda$-embedding from $B$ into $c_0$. Then there is a strict and good-$\lambda$-embedding from $M$ into $c_0$. 

\bn\rm
In particular, if $X$ is a Banach space and if there exists a good-$\lambda$-embedding
from its closed unit ball into $c_0$, then there exists a good-$\lambda$-embedding from $X$ 
into $c_0$. The following extension result is obvious.

\bn\bf
Proposition 5.2. \it Assume that $(M,d)$ is a separable metric space and that $N$ is a dense subset of $M$.
If there is a good-$\lambda$-embedding from $N$  into $c_0$, then there is a good-$\lambda$-embedding from $M$ into $c_0$.

\bn\bf
Remark 5.3. \rm In Definition 2.2, we didn't specify if the balls were closed or open.
We can define two different properties, $\pi(\lambda)$ with closed balls and $\pi(\lambda)$ with open balls.
These two properties are equivalent! 
Indeed, the proof of Proposition 2.3 shows that if there is a good-$\lambda$-embedding from $M$ into $c_0$,
then $M$ has property $\pi(\lambda)$ with closed balls, which in turn implies 
that $M$ has property $\pi(\lambda)$ with open balls.
On the othe hand, the proof of Theorem 4.1 shows that if $M$ has property $\pi(\lambda)$ with open balls, 
then  there is a good-$\lambda$-embedding from $M$ into $c_0$. This proves that 
property $\pi(\lambda)$ with open balls is equivalent to property $\pi(\lambda)$ with closed balls.

\bn\rm
N. Kalton and G. Lancien introduced the following definition :
A metric space $(M,d)$ has $\Pi(\lambda)$ if,
for every $\mu>\lambda$,
there exists $\nu>\mu$ such that  for every closed balls $B_1$ and $B_2$ with positive radii $r_1$ and $r_2$,
the exists subsets $U_1,\cdots U_N,V_1,\cdots,V_n$ of $M$ such that the sets
$U_n\times V_n$ are a covering of $E_\mu:=\{(x,y)\in B_1\times B_2;\,d(x,y)>\mu(r_1+r_2)\}$ and,
$$
\text{for all n, }\quad\lambda\delta(U_n,V_n)\ge\nu(r_1+r_2)
$$

\mn\bf
Lemma 5.4. \it Property $\Pi(\lambda)$ implies property $\pi(\lambda)$. 

\mn\rm
We do not know if the converse is true. Let us notice that N. Kalton and G. Lancien proved  that if a separable
metric space satisfies property $\Pi(\lambda)$, then there exists $f:M\to c_0$ such that, for all $x,y\in M$, $x\ne y$,
we have 
$$
d(x,y)<\Vert f(x)-f(y)\Vert\le\lambda d(x,y)
$$
which is  weaker than the condition $f$ is a strict and good-$\lambda$-embedding. 
Theorem 1 improve their result since our hypothesis, $M$ has $\pi(\lambda)$, is weaker,
and our conclusion, $f$ is a strict and good-$\lambda$-embedding, is stronger. Moreover, our condition $\pi(\lambda)$
is a necessary and sufficient condition for the existence of a good-$\lambda$-embedding.

\mn\sl Proof of Lemma 5.4. \rm Let us assume that $(M,d)$ has property $\Pi(\lambda)$.
Let $E\subset B_1\times B_2$ such that $\delta(E)>\lambda(r_1+r_2)$. We fix $\mu>\lambda$
such that $\delta(E)>\mu(r_1+r_2)$. Then $E\subset E_\mu$.  
Let $\nu>\mu$ be given by property $\Pi(\lambda)$. 
Let $1=a_1<a_2<\cdots<a_K$ be a sequence such that $diam(E_\mu)=a_K\mu(r_1+r_2)$ 
and $\frac{a_{k+1}}{a_{k}}<\frac{\nu}{\mu}$ whenever $1\le k< K$.
We denote 
$$
E_k:=\{(x,y)\in B_1\times B_2;\,a_{k}\mu(r_1+r_2)<d(x,y)\le \mu(r_1+r_2)a_{k+1}\}.
$$
The $E_k$'s form a covering of $E_\mu$.
Let $B_1^k$ and $B_
2^k$ be the closed balls of the same center as $B_1$ and $B_2$ and of radius
$a_kr_1$ and $a_kr_2$ respectively. Obviously, $E_k\subset B_1^k\times B_2^k$. 
Applying property $\Pi(\lambda)$ for each $k$,
We can find  subsets $U_{k,1},\cdots,U_{k,N_k}$, $V_{k,1},\cdots,V_{k,N_k}$ of $M$ such that
the sets $U_{k,n}\times V_{k,n}$ for $1\le n\le N_k$ form a covering of $E_k$ and
$$
\text{for all n, }\quad\lambda\delta(U_{k,n},V_{k,n})\ge\nu(a_kr_1+a_kr_2)>\mu a_{k+1}(r_1+r_2)
$$
We can assume in addition that for each $n$, the sets $U_{k,n}\times V_{k,n}$ are pairwise disjoint
(because a finite union of products can  always be written as a finite union of pairwise disjoint products).
If we denote $E_{k,n}=E\cap E_k\cap (U_{k,n}\times V_{k,n})$, the $E_{k,n}$'s form a partition of $E$.
Moreover,
$\pi(E_{k,n})\subset U_{k,n}\times V_{k,n}$, and the above inequality implies
$$
\text{for all n, }\quad\lambda\delta(\pi(E_{k,n}))>\mu a_{k+1}(r_1+r_2)\ge diam(E_k)\ge diam(E_{k,n}).
$$
We have proved property $\pi(\lambda)$.

\bn\mn\bf
6) Strict-$\lambda$-embeddings into $c_0^+$.

\bn\rm
Here, $c_0^+$ denotes the positive cone of $c_0$. 
Observe that :
$$
u,v\in c_0^+\quad\Rightarrow\Vert u-v\Vert\le\max\{\Vert u\Vert,\Vert v\Vert\}.
$$
The existence of a strict-$\lambda$-embedding into $c_0^+$
follows from the following property $\pi^+(\lambda)$.

\bn\bf
Definition 6.1. \rm A metric space $(M,d)$ has property  $\pi^+(\lambda)$ (with $\lambda>1$)
if,

\noindent
a) Whenever $B_1$ and $B_2$ are balls
of positive radii $r_1$ et $r_2$ and  $E$  is a subset of $B_1\times B_2$
such that $\delta(E)>\lambda\max(r_1,r_2)$,  
there exists a finite partition $\{E_1,\cdots,E_N\}$ of $E$ 
satisfying
$$
\text{ for each }n, \quad
diam(E_n)< \lambda\delta(\pi(E_n))
$$
b) There exists $\theta<\lambda$ and $\varphi:M\to[0,+\infty[$ such that \quad
$$
\vert\varphi(x)-\varphi(y)\vert\le d(x,y)\le\theta\max\bigl(\varphi(x),\varphi(y)\bigr) \quad\text{for all }x,y\in M.
$$
The function $\varphi$ is called a control function.

\bn\bf
Remark 6.2. \rm
1) It is easy to see that it is enough to check a) whenever $r_1=r_2(=r)$.

\sn
2) If $\lambda>2$, the function $\varphi(x)=d(x,a)$ 
is a control function (take $\theta=2$). Therefore, the metric space $M$ has property $\pi^+(\lambda)$ 
if and only if the bounded subsets of $M$ have property $\pi^+(\lambda)$.
In particular a Banach space $X$ has property $\pi^+(\lambda)$ (with $\lambda>2$)
if and only if its unit ball  has property $\pi^+(\lambda)$.

\sn
3) If $M$ is bounded, then, for any $\lambda>1$, the function
$\varphi:M\to[0,+\infty[$ given by $\varphi(x)=d(x,a)+diam(M)$
satisfies condition b) of property $\pi^+(\lambda)$ (take $\theta=1$).

\bn\bf
Proposition 6.3. \it
1) If there is a $\lambda$-embedding from $(M,d)$ into $c_0^+$, then
$M$ has property 
$\pi^+(\mu)$ for all $\mu>\lambda$.

\sn
2) Assume that $M$ is bounded or that $\lambda>2$. If there is a good-$\lambda$-embedding from $M$ into $c_0$,
then $M$ has property~$\pi(\lambda)$.

\mn\sl Proof. \rm
Let $B_1=B(a_1,r)$ and $B_2=B(a_2,r)$. Let $E\subset B_1\times B_2$
such that  $\lambda r+\eps<\delta(E)$ for some $\eps>0$. 
Let $f:M\to c_0$ be a $\lambda$-embedding given by $f(x)=\sum\limits_{i=0}^{+\infty}f_i(x)e_i$. We denote $P(x)=\sum\limits_{i=0}^{i_0}f_i(x)e_i$ and $Q=f-P$,
where $i_0$ is such that $\max\{\Vert Q(a_1)\Vert,\Vert Q(a_2)\Vert\}<\eps$.
If $(x,y)\in E$, then
$$
\eqalign{\Vert Q(x)-Q(y)\Vert&\le\max\{\Vert Q(x)\Vert,
\Vert Q(y)\Vert\}\cr
&\le\max\big\{\Vert Q(x)-Q(a_1)\Vert+\Vert Q(a_1)\Vert,\Vert Q(y)-Q(a_2)\Vert+\Vert Q(a_2)\Vert\big\}\cr
&\le\lambda r+\eps<\delta(E)\le d(x,y)\le\Vert f(x)-f(y)\Vert.
\cr}
$$
Thus,
$\Vert f(x)-f(y)\Vert=\Vert P(x)-P(y)\Vert$. Following the lines of the proof of Proposition~2.3 we get,
for any $\mu>\lambda$, a partition $\{E_1,\cdots,E_N\}$ of $E$ such that 
for each $n$, $diam(E_n)\le\mu\delta(\pi(E_n))$.
Condition a) can be checked and the $E_n$'s satisfy also 
$\delta(\pi(E_n))>r$. Moreover, if $\varphi(x)=\Vert f(x)\Vert/\lambda$, then $\vert\varphi(x)-\varphi(y)\vert\le d(x,y)\le\lambda\max\bigl(\varphi(x),\varphi(y)\bigr)$ 
\quad for all $x,y\in M$. This proves that $\varphi$ is a control function of $\pi^+(\mu)$ because $\lambda<\mu$.

\noindent The proof of 2) also follows the lines of the corresponding case of Proposition~2.3,
and here we do not have to worry about the existence of a control function by Remark 6.2.

\bn\bf
Corollary  6.4. \it Let $X$ be a Banach space. If there exists $u\in S_X$
and an infinite dimensional subspace $Y$ of $X$ such that
$\inf\{\Vert u+2y\Vert;\,y\in S_Y\}>\lambda$,
then there is no $\lambda$-embedding from $M$ into $c_0^+$.

\mn\sl Proof. \rm
If $E=\{(u+2y,-u-2y);\,y\in S_Y\}\subset\overline{B}(u,2)\times\overline{B}(-u,2)$, 
$E$ satisfies $\delta(E)>2\lambda$.
If there is a $\lambda$-embedding from $X$ into $c_0^+$, 
then, by Proposition 6.3, there is a partition  
$\{E_1,\cdots,E_N\}$ of $E$  such that for each  $n$,\quad
$\delta(\pi(E_n))>2$. If $F_n$ is the norm closure of $E_n$, then 
$\{F_1,\cdots,F_N\}$ is a covering of $E$ and we still have $\delta(\pi(F_n))>2$.

\sn
But as in the proof of Corollary 2.4, we also have $\delta(\pi(F_n))\le 2$,
which is absurd.

\bn\bf Example 6.5. \it The metric space $\ell^1$ 
does not  $\lambda$-embed into $c_0^+$ whenever 
$\lambda<3$. \rm (If $u=e_0$ and $Y=\{y=(y_i)\in\ell^1;\,y_0=0\}$, then 
$\Vert u+2y\Vert=3$ for each $y\in S_Y$).

\sn\it
The space $\ell^p$ does not $\lambda$-embed into $c_0^+$ whenever 
$\lambda<(1+2^p)^{1/p}$ (since in this case, for all $y\in S_Y$, $\Vert u+2y\Vert=(1+2^p)^{1/p}$).

\sn
The balls of positive radius of $\ell^p$ do not $\lambda$-embed into $c_0^+$ whenever 
$\lambda<(1+2^p)^{1/p}$. On the other hand, it follows from Corollary 4.4 that the balls of 
positive radius of $\ell^p$ embed in $c_0$ if $\lambda=2^{1/p}<(1+2^p)^{1/p}$.

\sn\rm
Indeed, assume that for some $\lambda<(1+2^p)^{1/p}$, a ball $B$ 
of positive radius of $\ell^p$ $\lambda$-embeds into $c_0^+$.
We can assume that $B$ is the unit ball of $\ell^p$. According to Proposition 6.3,
$B$ has property $\pi^+(\mu)$ for every $\mu>\lambda$, and by Remark 6.2 2), $\ell^p$
has property $\pi^+(\mu)$ for every $\mu>\lambda$, and from Theorem 6.9 below,
$\ell^p$ $\mu$-embeds into $c_0^+$ for every $\mu>\lambda$. But this contradicts the fact 
that $\ell^p$ does not $\mu$-embed into $c_0^+$ whenever $\mu<(1+2^p)^{1/p}$.

\bn\bf
Example 6.6. \it A compact metric space $M$ has property $\pi^+(\lambda)$ for all $\lambda>1$.
A metric space $M$ such that its bounded subsets are totally bounded has property $\pi^+(\lambda)$ for all $\lambda>2$,
but may fail property $\pi^+(2)$.

\mn\rm
A metric space $M$ such that its bounded subsets are totally bounded satisfies condition a) of property $\pi^+(\lambda)$ 
for all $\lambda>1$,
and any metric space satisfies condition b) of property $\pi^+(\lambda)$ for all $\lambda>2$.
Moreover, if $\lambda>1$ and $M$ is compact, 
then $M$ is bounded and so satisfies condition b). 

\sn
The bounded subsets of the set \bf Z \rm of integers are finite. 
Let $\varphi:$\bf Z\rm$\to[0,+\infty[$ such that
$\vert\varphi(x)-\varphi(y)\vert\le \vert x-y\vert$ for all $x,y\in$\bf Z\rm. 
If $x_n=(-1)^nn$, then 
$\vert x_n-x_{n+1}\vert=2n+1$ and $\varphi(x_n)\le\varphi(0)+n$. 
Consequently, if $\theta<2$, then   
$\vert x_n-x_{n+1}\vert>\theta\max\big\{\varphi(x_n),\varphi(x_{n+1})\big\}$ for $n$ large enough. Thus \bf Z \rm
do not admit any control function $\varphi$
for property $\pi^+(2)$.

\bn\bf
Example 6.7. \it Each metric space $M$ has property $\pi^+(3)$.

\mn\sl Proof. \rm
Let $B_1$ and $B_2$ 
be balls of radius~$r$ and   $E\subset B_1\times B_2$
such that $\eps:=\delta(E)-3r>0$.
As in Example $2$, using the fact, there is a partition $\{E_1,\cdots,E_N\}$
of $E$ such that for each $n$,
 $\delta(E_n)+\eps>diam(E_{n})$.
If $(u,v)\in \pi(E_n)$, there is $v'\in B_2$
so that $(u,v')\in E_n$. Moreover $v,v'\in B_2$, hence~:  
$$
3d(u,v)\ge 3d(u,v')-3d(v',v)\ge 3\delta(E_n)-6r\ge\delta(E_n)+2\eps> diam(E_n)
$$
Taking the infimum over all $(u,v)\in \pi(E_n)$, we get
$3\delta(\pi(E_n))\!>\! diam(E_{n})$.

\bn\bf
Example 6.8. \it $\ell^p(\entier)$ has property $\pi^+\bigl((1+2^p)^{1/p}\bigr)$.

\mn\sl Proof. \rm
Let   
$E\subset B(a_1,r)\times B(a_2,r)$ such that $\alpha=\delta(E)^p-(1+2^p)r^p>0$. We choose  $\eps>0$
such that $(2r+\eps)^p<\frac{2^p}{1+2^p}\delta(E)^p-\alpha/2$ and $(t+\eps)^p-\alpha/2<t^p$ si $0\le t\le diam(E)$. 
Let $(e_i)$ be the canonical basis of $\ell^p$, $P,Q:\ell^p\to\ell^p$
be defined by $P(\sum\limits_{i=0}^\infty x_ie_i)=\sum\limits_{i=0}^{i_0} x_ie_i$
and $Q=I-P$,
where $i_0$ is chosen so that $\Vert Qa_1-Qa_2\Vert<\eps$. 
Since $Lip(P)=1$ and $P$ has its values in a finite dimensional space,
Fact 2.1 implies the existence of a partition $\{E_1,\cdots,E_N\}$ of $E$ such that, 
for each $n$, if $(x,y)\in E_n$, then
$\Vert Px-Py\Vert\le\delta(\pi(E_n))+\eps$. 
On the other hand, 
$\Vert Qx -Qy\Vert\le\Vert Qx-Qa_1\Vert+\Vert Qy-Qa_2\Vert+\Vert Qa_1-Qa_2\Vert\le 2r+\eps$.
Hence
$$
diam(E_n)^p\le\bigl(\delta(\pi(E_n))+\eps\bigr) ^p+(2r+\eps)^p\le\bigl(\delta(\pi(E_n))+\eps\bigr) ^p+ \frac{2^p}{1+2^p}\delta(E)^p-\alpha/2
$$
and so $diam(E_n)^p\le(1+2^p)\bigl(\delta(\pi(E_n))+\eps\bigr) ^p-(1+2^p)\alpha/2<(1+2^p)\delta(\pi(E_n))^p$.

\bn\bf
Theorem 6.9. \it If the separable metric space $(M,d)$  has property $\pi^+(\lambda)$ with $1\!<\!\lambda\!\le\! 3$, 
then there exists $f:M\to c_0^+$ such that
for all $x,y\in M$, $x\ne y$, we have~:
$$
d(x,y)<\Vert f(x)-f(y)\Vert<\lambda d(x,y).
$$

\mn\bf
Corollary 6.10. \it Let $(M,d)$ be a separable metric space. 
Then there exists $f:M\to c_0^+$ such that,
for all $x,y\in M$, $x\ne y$, we have~:
$$
d(x,y)<\Vert f(x)-f(y)\Vert<3 d(x,y).
$$

\rm\noindent
This result is optimal since we observed that there is no $\lambda$-embedding from $\ell^1$
into $c_0^+$ with $\lambda<3$.

 \bn\bf
Corollary 6.11. \it 
There exists $f:\ell^p\to c_0^+$ such that, for all $x,y\in M$, $x\ne y$, we have~:
$$
d(x,y)<\Vert f(x)-f(y)\Vert<(1+2^p)^{1/p} d(x,y).
$$

\rm\mn
This result is optimal since we observed that there is no $\lambda$-embedding from $\ell^p$
into $c_0^+$ whenever $\lambda<(1+2^p)^{1/p}$.

\bn\bf
Corollary  6.12. \it  If $(M,d)$ is a compact space and $\eps>0$,
there exists $f:M\to c_0^+$ such that, for all $x,y\in M$, $x\ne y$, we have~:
$d(x,y)<\Vert f(x)-f(y)\Vert<(1+\eps)d(x,y).$

\mn
If $(M,d)$ is a metric space such that its bounded subsets are totally bounded,
there exists $f:M\to c_0^+$ such that, for all $x,y\in M$, $x\ne y$, we have~:
$d(x,y)<\Vert f(x)-f(y)\Vert<(2+\eps) d(x,y).$

\rm\mn
The following result shows that  we cannot replace $2+\eps$ by $2$ in the above statement.

\bn\bf
Proposition 6.13. \it There exists a separable metric space $M$
such that, for any $\lambda>1$, $M$ $\lambda$-embeds into $c_0$ but 
there is no $2$-embedding from $M$ into $c_0^+$.

\mn\sl
Proof. \rm
Let $(e_n)$ be the canonical basis of $\ell^1(\entier)$ 
and $F_p:=\{pe_k,\,e_0+pe_k;\,1\le k\le p\}$.
We define $M=\{0,e_0\}\cup\bigcup\limits_{p=1}^{+\infty}F_p\,\subset\,\ell^1(\entier)$.
The bounded sets of $M$ are finite, hence totally bounded, so, by Corollary 2,
for any $\lambda>1$, $M$ $\lambda$-embeds into $c_0$.

\sn
Assume now that there exists $f=(f_n):M\to c_0^+$ such that, for all $x,y\in M$,
$$
\Vert x-y\Vert_1\le\Vert f(x)-f(y)\Vert_\infty\le 2\Vert x-y\Vert_1
$$
Let us denote $C=\max\{\Vert f(0)\Vert_\infty,\Vert f(e_0)\Vert_\infty\}$, fix $n_0\ge 1$
such that for all $n>n_0$, one has $f_n(0)<1$ and $f_n(e_0)<1$,
and finally fix $p>C/2+1$.
We claim that the mapping $\varphi:\{1,\cdots,p\}\to\{0,1\}^{n_0}$
defined by $\varphi(k)=\bigl(\un_{[0,C+1]}(f_n(pe_k))\bigr)_{n\le n_0}$ is injective.
This leads to a contradiction if we also have $p>2^{n_0}$.

\sn
If $n\ge 0$ and $1\le k\le p$, then  
$$
f_n(pe_k)\le\vert f_n(pe_k)-f_n(0)\vert+f_n(0)
\quad\text{and}\quad
f_n(pe_k+e_0)\le\vert f_n(pe_k+e_0)-f_n(e_0)\vert+f_n(e_0)
$$
so 
$$
f_n(pe_k)\le 2p+C\qquad\text{and}\qquad f_n(e_0+pe_k)\le 2p+C.\tag 1
$$ 
Whenever $n>n_0$, we have a better estimate:
$$
f_n(pe_k)<2p+1\qquad\text{and}\qquad f_n(e_0+pe_k)<2p+1
$$
So, if $n>n_0$, 
$$
\vert f_n(pe_k)- f_n(e_0+pe_k)\vert\le\max\{f_n(pe_k),\, f_n(e_0+pe_k)\}<2p+1
$$
On the other hand, if $1\le k\ne\ell\le p$, we have 
$$
2p+1=\Vert e_0+pe_k-pe_\ell\Vert_1\le\Vert f(e_0+pe_k)-f(pe_\ell)\Vert_\infty,
$$
hence, there exists $n\le n_0$ such that $\vert f_n(e_0+pe_k)-f_n(pe_\ell)\vert\ge 2p+1$,
and using the fact that $\vert f_n(e_0+pe_k)-f_n(pe_k)\vert\le 2$, we obtain 
$$
\vert f_n(pe_k)-f_n(pe_\ell)\vert\ge 2p-1\tag 2
$$
Using (1) and (2), we obtain that either $f_n(pe_k)\le C+1$ and $f_n(pe_\ell)\ge 2p-1$,
or $f_n(pe_\ell)\le C+1$ and $f_n(pe_k)\ge 2p-1$,
hence $\un_{[0,C+1]}(f_n(pe_k))\ne\un_{[0,C+1]}(f_n(pe_\ell))$, and $\varphi$ is injective.

\bn\rm
The proof of Theorem 6.9 is analogous to the proof of Theorem 4.1
and relies on the following two lemmas (analogous to Lemmas 4.5 and 4.6).

\bn\bf
Lemma 6.14. \it Let $(M,d)$ be a metric space, $U,V,F$ non empty bounded subsets of $M$ and $\eps\ge 0$.
There exists $f:M\to\reel^+$, such that $Lip(f)\le 1$ and~:

\medskip
1) For all $x\in F$, \quad $f(x)\le\eps$,

\medskip
2) For all $(x,y)\in U\times V$, \quad $f(x)-f(y)=\min\big\{\delta(U,V),\max(\delta(U,F),\delta(V,F))+\eps\big\}$.

\mn\sl Proof. \rm
Indeed, if $\delta(V,F)\le\delta(U,F)$ and if we put  $t=\min(\delta(U,V),\delta(U,F)+\eps)$, the function $f$ defined by  
$f(x)=\max(t-d(x,U),0)$ satisfies Lemma 6.14.

\bn\bf
Lemma 6.15. \it Let $(M,d)$ be a metric space with property $\pi^+(\lambda)$ with $1<\lambda\le 3$, 
$F\subset G$ be finite subsets of $M$ and $0<\alpha<\beta$, we set~:
$$
A_+(G,\alpha)=\big\{(x,y)\in M\times M;\,
d(x,y)\ge\lambda\bigl(\max(d(x,G),d(y,G))+\alpha\bigr)\big\}
$$
Then there exists a finite partition $\{E_1,\cdots,E_N\}$ of $A_+(G,\alpha)\backslash A_+(F,\beta)$
such that,  if we denote $\pi(E_n)=U_n\times V_n$,
then
$$
\text{for each }n,\qquad diam(E_n)<\lambda\min\big\{\delta(U_n,V_n),\max(\delta(U_{n},F),\delta(V_{n},F))+2\beta\big\}.
$$

\mn\sl Proof. \rm
Denote $K=\sup\{\vert\varphi(a)\vert;\,a\in G\}$, where $\varphi$ is the control function. 
For all $x\in M$, we have $\varphi(x)\le d(x,G)+K$. 
If $(x,y)\in A_+(G,\alpha)$, then
$$
\eqalign{
\lambda\max\bigl(d(x,G),d(y,G)\bigr)&\le d(x,y)\le\theta\max\bigl(\varphi(x),\varphi(y)\bigr)\cr
&\le\theta\max\bigl(d(x,G),d(y,G)\bigr)+\theta K.\cr
}
$$
Since $G$ is bounded and $\lambda>\theta$, we can find  $B\subset M$ bounded such that 
$A_+(G,\alpha)\subset B\times B$. 
The rest of the proof follows the lines of  Lemma 4.6 (show that  $\delta(E_{jk})>\max\{r_j,r_k\}$).

\mn 
For the proof of Theorem 6.9, choose $\eps_1\!>\!\varphi(a_1)$, 
which implies, for all $x,y\!\in\! M$, $d(x,y)\le\lambda\max(\varphi(x),\varphi(y))<\lambda\bigl(\max(d(x,a_1),d(y,a_1))+\eps_1\bigr):=\sigma_1$.

\bn\bf
Remark 6.16. \rm Property $\pi(\lambda)$ characterizes the existence of a good-$\lambda$-embedding
into $c_0$, but we do not know if property $\pi^+(\lambda)$ characterizes the existence of a good-$\lambda$-embedding
into $c_0^+$. We do not know of any internal characterization of the existence of a good-$\lambda$-embedding
into $c_0^+$, or of the existence of a $\lambda$-embedding
into $c_0^+$. However, it seems very likely that if $M$ is bounded, or if $2<\lambda\le 3$, then
the existence of a good-$\lambda$-embedding into $c_0^+$ is equivalent to the fact that $M$ has property $\pi^+(\lambda)$.

\bn\it
We wish to thank anonymous referees for their useful comments on the presentation of this paper.

\bn\bn\bf
References.

\mn\rm
[1] I. Aharoni, Every separable metric space is Lipschitz equivalent to a subset of 
$c_0^+$, Israel J. Math. 19 (1974), 284-291.

\mn
[2] P. Assouad, Remarques sur un article de Israel Aharoni sur les prolongements lipschitziens dans $c_0$, 
Israel J. Math. 31 (1978), 97-100.

\mn
[3] N. J. Kalton and G. Lancien, Best constants for Lipschitz embeddings of metric spaces
into $c_0$,  Fundamenta Mathematicae, 199, 2008, 249-272.

\mn
[4] J. Pelant, Embeddings into $c_0$, Topology Appl. 57 (1994), 259-269.

\enddocument